\documentclass{article}
\usepackage{maa-monthly}

\usepackage{xcolor}

\theoremstyle{theorem}
\newtheorem{theorem}{Theorem}
 
\theoremstyle{definition}
\newtheorem{definition}[theorem]{Definition}
\newtheorem{remark}[theorem]{Remark}

\newcommand{\fs}{\mathop{\mathrm{FS}}}
\newcommand{\fr}{\mathnormal{\mathrm{Fr}}}
\newcommand{\zf}{\mathnormal{\mathsf{ZF}}}

\definecolor{rojo}{HTML}{cc0000}
\definecolor{azul}{HTML}{000099}

\renewcommand{\r}{{\color{rojo}\mathnormal{\mathrm{red}}}}
\renewcommand{\b}{{\color{azul}\mathnormal{\mathrm{blue}}}}
\newcommand{\tr}{{\color{rojo}red}}
\newcommand{\tb}{{\color{azul}blue}}
\DeclareMathOperator{\ev}{ev}

\theoremstyle{theorem}
\newtheorem{proposition}[theorem]{Proposition}
\newtheorem{lemma}[theorem]{Lemma}
\newtheorem{corollary}[theorem]{Corollary}

\begin{document}

\title{Using Ultrafilters to Prove \\ Ramsey-type Theorems}
\markright{Ultrafilters and Ramsey theory}
\author{David J. Fern\'andez-Bret\'on}

    \begingroup
    \renewcommand{\thefootnote}{}
    \endgroup

\maketitle

\begin{abstract}
Ultrafilters are a tool, originating in mathematical logic and general topology, that has steadily found more and more uses in multiple areas of mathematics, such as combinatorics, dynamics, and algebra, among others. The purpose of this article is to introduce ultrafilters in a friendly manner and present some applications to the branch of combinatorics known as Ramsey theory, culminating with a new ultrafilter-based proof of van der Waerden's theorem.
\end{abstract}

\section{Introduction.}

\noindent Learning the extent to which ultrafilters are useful in multiple branches of mathematics can be rather surprising, given the scarcity of occasions when ultrafilters show up throughout the usual undergraduate and graduate curriculum. A typical undergraduate student will for the most part finish her degree without ever encountering the definition of an ultrafilter, and those few who do will probably do so in the context of general topology, where they will very likely be told that ultrafilters are one possible way of generalizing the notion of a convergent sequence. Slightly more esoteric is the presence of ultrafilters in courses on model theory or set theory, which are subjects that relatively few undergraduates, or even graduate students, ever undertake for study. It has slowly become more apparent, however, that ultrafilters are an extremely useful tool for a variety of different purposes in mathematics. Nowadays, there is abundant literature about the interaction between ultrafilters and combinatorics, particularly additive combinatorics (see, e.g.,~\cite{hindman-number-theory,hindman-update,hindman-strauss}). Also classical are the applications of ultrafilters to algebraic geometry, which for the most part occur in the context of model theory (especially applications of ultraproducts---which are a construction that arises naturally from ultrafilters---such as the Ax--Grothendieck theorem~\cite{ax-finite-fields} and the Ax--Kochen theorem~\cite{ax-kochen}). More recently, new applications have arisen in commutative algebra, for example, Schoutens's work using ultraproducts of commutative rings for a variety of purposes, e.g., generalizing the notion of tight closure to certain rings of characteristic 0~\cite{schoutens}. And it is impossible not to mention also the many recent applications of methods from nonstandard analysis---a theory which is tightly connected with ultrafilters---in a variety of combinatorial settings~\cite{dinasso-goldbring-lupini}.

The main purpose of this article is to introduce a hypothetical reader who has never encountered ultrafilters to the basic notions about them, mentioning enough results so that some applications can be properly appreciated. In regard to the theory of ultrafilters, we provide proofs for only a few select statements, those that are sufficiently simple to prove and useful for developing intuition. The remaining statements will only be mentioned without proofs; on the other hand, with regard to the applications of these tools to Ramsey theory, we provide full proofs in complete detail. This way we hope that the reader will, after reading this article, be able to both understand other proofs of Ramsey-theoretic statements that use ultrafilters, and have the tools to attempt their own such proofs.

This article is organized in an oscillatory way, in the sense that we continuously alternate between stating definitions and standard facts about ultrafilters, and providing a paradigmatic application of the newly introduced ideas. In Section 2, we begin by introducing the basic definitions concerning ultrafilters, as viewed combinatorially; after a brief digression (a subsection where we prove the existence of nonprincipal ultrafilters), we show how this combinatorial device can contribute to proving Ramsey's theorem. In Section 3, we start by introducing some more theory, focused mostly on the definition of a right-topological semigroup structure on the space of all ultrafilters; immediately after that we follow up by showing (with full detail) how this theory can be applied to prove Schur's theorem, and hint at how the same ideas can be used to prove the Folkman--Rado--Sanders theorem and Hindman's theorem, which are generalizations of Schur's theorem. Finally, in Section 4, at the beginning we introduce yet more theory, particularly about minimal idempotents and minimal ideals in the semigroup of ultrafilters, in order to immediately follow up with an example of how these newly introduced concepts can be applied, concretely, to prove van der Waerden's theorem. The proofs of the applications from Sections 1 and 2 are standard (although we believe that the way we sketch the proof of Hindman's theorem differs slightly from other presentations currently available); meanwhile, the specific proof of van der Waerden's theorem presented here---different from other ultrafilter proofs of the same result, e.g., the one in~\cite[Theorem 14.1 \& Corollary 14.2]{hindman-strauss}---is, to the best of our knowledge, completely new.

\section{Basic definitions, colorings of graphs.}

In this section we introduce the definition of ultrafilters, attempting to provide some intuition along the way, and then show how to apply these concepts to provide neat proofs of various versions of Ramsey's theorem, which deals with colorings of graphs. We also provide a self-contained proof that ultrafilters exist, but the reader who is mostly interested in applications may skip the corresponding subsection without harm.

\subsection{Definition and basic properties of ultrafilters.}

\begin{definition}\label{def:ultrafilter}
Given a set $X$, an \textbf{ultrafilter} on $X$ is a family $u\subseteq\wp(X)$ satisfying 
\begin{enumerate}
\item For every $A,B\subseteq X$, $A\cap B\in u$ if and only if $A\in u\text{ and }B\in u$; 
\item for every $A,B\subseteq X$, $A\cup B\in u$ if and only if $A\in u\text{ or }B\in u$;
\item for every $A\subseteq X$, $A\in u$ if and only if $X\setminus A\notin u$.
\end{enumerate}
\end{definition}

From the perspective of mathematical logic, the intuitive idea of an ultrafilter is that it is a way of assigning truth values to all subsets of $X$ (here, we interpret $A\in u$ as ``$A$ is true'' and $A\notin u$ as ``$A$ is false''). From this perspective, the definition simply says that ultrafilters respect the fact that conjunction corresponds to intersection, disjunction corresponds to union, and negation corresponds to complement.\footnote{As a matter of fact, if one sees $\wp(X)$ as a Boolean algebra, equipped with the operations of union for joins, intersection for meets, and complement for negatives, then one can clearly see that the definition stated above is equivalent to defining an ultrafilter as the preimage of $1$ under a Boolean algebra homomorphism $\wp(X)\longrightarrow\{0,1\}$.} If one prefers to think in terms of measure theory, then one can conceive of an ultrafilter as a finitely additive $\{0,1\}$-valued measure on $\wp(X)$ (the power set of $X$), the translation being given by interpreting the measure of a subset $A\subseteq X$ as $1$ if and only if $A\in u$, and $0$ otherwise.

The logical interpretation makes it clear that, in Definition~\ref{def:ultrafilter}, one could assume only conditions (1) and (3), and this would automatically imply (2)---alternatively, one could assume (2) and (3) and this would imply condition (1)---because of the well-known classical result in propositional logic that every binary connective can be written in terms of conjunction and negation; alternatively, every binary connective can be written in terms of disjunction and negation. In fact, every binary connective can be written in terms of the Sheffer stroke, which means that we could have equivalently defined an ultrafilter as a family $u$ of subsets of $X$ satisfying that
\begin{equation*}
X\setminus(A\cap B)\in u\text{ if and only if }(\text{not both }A\in u\text{ and }B\in u),
\end{equation*}
for all $A,B\subseteq X$.

Notice that, if $X=\varnothing$, then there are no ultrafilters over $X$ (for the definition would require that, if $u$ is an ultrafilter on $X=\varnothing$, then $\varnothing\in u$ if and only if $\varnothing=X\setminus\varnothing\notin u$). Thus, from now on, we focus on ultrafilters over nonempty sets. Observe that if $X\neq\varnothing$ and $u$ is an ultrafilter on $X$ then $u$ is nonempty (because now $\varnothing$ and $X=X\setminus\varnothing$ are distinct, and by condition (3) in Definition~\ref{def:ultrafilter} we must have one of these two sets belonging to $u$).

\begin{proposition}\label{prop:closedupwards}
Every ultrafilter is closed upwards. That is, if $X$ is a set and $u$ is an ultrafilter on $X$, then whenever $A\in u$ and $A\subseteq B\subseteq X$ we must have $B\in u$. In particular, we always have $X\in u$ (and consequently $\varnothing\notin u$).
\end{proposition}

\begin{proof}
Suppose that $A\in u$ and $A\subseteq B\subseteq X$. Since $A\in u$, we (quite tautologically) have that $A\in u$ or $B\in u$. Hence, by condition (2) in Definition~\ref{def:ultrafilter}, it must be the case that $B=A\cup B\in u$.

The claim in the second sentence of Proposition~\ref{prop:closedupwards} follows immediately from the fact that every ultrafilter is nonempty, as remarked right before the statement of this proposition.
\end{proof}

\begin{proposition}\label{prop:finiteunion}
Let $X$ be a set, and let $u$ be an ultrafilter on $X$. If $A\in u$ and $A=A_1\cup\cdots\cup A_n$, then there exists an $i$ such that $A_i\in u$; moreover, if the $A_j$ are pairwise disjoint then this $i$ is unique. In particular, for every finite partition of $X$, exactly one piece of the partition belongs to $u$.
\end{proposition}

\begin{proof}
The first part of the proposition is a straightforward induction, where the base case $n=2$ follows directly from condition (2) in Definition~\ref{def:ultrafilter}. For the second part, notice that if the $A_j$ were pairwise disjoint and we had $A_i,A_k\in u$ for $i\neq k$, then we would also have $\varnothing=A_i\cap A_k\in u$, contradicting Proposition~\ref{prop:closedupwards}.
\end{proof}

The easiest example of an ultrafilter is the following: if $X$ is a set and $x\in X$, then the family $\{A\subseteq X\mid x\in A\}$ is an ultrafilter on $X$, as the reader should be able to verify straightforwardly. Ultrafilters fitting this description receive a special name.

\begin{definition}
Given a set $X$, an ultrafilter $u$ on $X$ is said to be {\bf principal} if there exists an $x\in X$ with $u_x:=\{A\subseteq X\mid x\in A\}=u$. Otherwise, $u$ will be said to be {\bf nonprincipal}.
\end{definition}

If one thinks in measure-theoretic terms, principal ultrafilters would correspond to {\em Dirac measures}: measures that are concentrated on a single point $x$, in the sense that the measure of any set $A$ is 1 if and only if $x\in A$, and 0 otherwise. The next theorem precisely characterizes principal ultrafilters.

\begin{theorem}\label{thm:principaluf}
Let $X$ be a set and let $u$ be an ultrafilter on $X$. Then $u$ is principal if and only if there exists a finite $F\subseteq X$ with $F\in u$.
\end{theorem}

\begin{proof}
If $u$ is principal, then there is an $x\in X$ such that $u=u_x=\{A\subseteq X\mid A\in u\}$; in particular $\{x\}\in u_x=u$. Conversely, suppose that $F\subseteq X$ is a finite set (say $F=\{x_1,\ldots,x_n\}$) with $F\in u$. Then $F=\{x_1\}\cup\cdots\cup\{x_n\}$, so by Proposition~\ref{prop:finiteunion}, there is an $i\in\{1,\ldots,n\}$ such that $\{x_i\}\in u$. We now claim that $u=u_{x_i}$: to see this, let $A\subseteq X$ and note that $A\in u$ if and only if $A\in u$ and $\{x_i\}\in u$, which in turn happens (by condition (1) in Definition~\ref{def:ultrafilter}) if and only if $A\cap\{x_i\}\in u$. Notice that 
\begin{equation*}
A\cap\{x_i\}=\begin{cases}\{x_i\}\text{ if }x_i\in A \\ \varnothing\text{ otherwise.}\end{cases}
\end{equation*}
Since $\varnothing\notin u$, it follows that $A\in u$ if and only if $x_i\in A$, so that $u=u_ {x_i}$.
\end{proof}

This immediately raises the question whether there is a set $X$ over which a nonprincipal ultrafilter exists. Such a set must be infinite by Theorem~\ref{thm:principaluf}. We devote the next subsection to proving that, for every infinite set $X$, there exists a nonprincipal ultrafilter $u$ over $X$. The reader who wishes to take this statement on faith and continue with applications is advised to skip the subsection.

\subsection{Existence of nonprincipal ultrafilters.}

Throughout this subsection, $X$ is assumed to be an infinite set. The objective of the section is to prove the existence of a nonprincipal ultrafilter on $X$.

\begin{definition}
We will say that a family $\mathcal F$ of subsets of $X$ is {\bf nice} if it is nonempty, closed under intersections (i.e., if $A,B\in\mathcal F$ then $A\cap B\in\mathcal F$), and $\varnothing\notin\mathcal F$.
\end{definition}

Notice that the collection of all nice families is partially ordered by inclusion. It can be readily checked that this collection thus partially ordered satisfies the hypothesis of Zorn's lemma (that every chain has an upper bound, in this case given by the union of the chain), and hence there exist maximal nice families. Moreover, whenever $\mathcal F$ is a nice family, there exists a maximal nice family $\mathcal M$ with $\mathcal F\subseteq\mathcal M$. We will use this fact to prove the existence of nonprincipal ultrafilters over $X$, which will follow immediately after the following two lemmas.

\begin{lemma}\label{lem:annexion}
Let $\mathcal M$ be a maximal nice family. If a set $A\subseteq X$ intersects every element of $\mathcal M$, then $A\in\mathcal M$. In particular, $M$ is closed upwards.
\end{lemma}

\begin{proof}
Since $\{X\}\cup\mathcal M$ is also a nice family containing $\mathcal M$, by maximality we have $\mathcal M=\mathcal M\cup\{X\}$ and so $X\in\mathcal M$. Now suppose that $A\subseteq X$ intersects every $B\in\mathcal M$. Then the family
\begin{equation*}
\{A\cap B\mid B\in\mathcal M\}\cup\mathcal M
\end{equation*}
is nice and contains $\mathcal M$, so by maximality $\mathcal M=\mathcal M\cup\{A\cap B\mid B\in\mathcal M\}$ and, in particular (since $X\in\mathcal M$), $A=A\cap X\in\mathcal M$. This finishes the proof of the first statement. For the second statement, if $A\in\mathcal M$ and $A\subseteq B$, then for every $C\in\mathcal M$ we have $B\cap C\supseteq A\cap C\in\mathcal M$; since $\varnothing\notin\mathcal M$ we get $B\cap C\neq\varnothing$. By the first statement, it follows that $B\in\mathcal M$.
\end{proof}

\begin{lemma}\label{lem:maximalnice}
A family of subsets of $X$ is maximal nice if and only if it is an ultrafilter.
\end{lemma}

\begin{proof}
We begin by proving the reverse implication. If $u$ is an ultrafilter, it follows immediately that $u$ is a nice family. Now, suppose that $\mathcal F$ is another nice family with $u\subseteq\mathcal F$. If the inclusion were proper, taking any $A\in\mathcal F\setminus u$ we would have that $X\setminus A\in u\subseteq\mathcal F$, and therefore $\varnothing=A\cap(X\setminus A)\in\mathcal F$, a contradiction. Hence, $u=\mathcal F$ and we are done.

Now for the forward implication, let $\mathcal M$ be a maximal nice family. 
We shall use the Sheffer stroke characterization of ultrafilters (see the second paragraph after Definition~\ref{def:ultrafilter}), that is, we will prove that $\mathcal M$ is an ultrafilter by showing that, for $A,B\subseteq X$, we have $X\setminus(A\cap B)\in\mathcal M$ if and only if not both of $A\in\mathcal M$, $B\in\mathcal M$ hold. So suppose that $X\setminus(A\cap B)\in\mathcal M$. Since $\varnothing\notin\mathcal M$ and $\mathcal M$ is closed under intersections, this means we must have $A\cap B\notin\mathcal M$. By Lemma~\ref{lem:annexion} this implies that there exists a $C\in\mathcal M$ that is disjoint from $A\cap B$. If $A\in\mathcal M$ then $A\cap C\in\mathcal M$ is disjoint from $B$; since $\mathcal M$ is closed under intersections and does not contain $\varnothing$, we must have that $B\notin\mathcal M$. Following the same reasoning with the roles of $A$ and $B$ interchanged shows that if $B\in\mathcal M$ then $A\notin\mathcal M$; in either case we have proved that not both of $A\in\mathcal M$, $B\in\mathcal M$ hold. Conversely, suppose that not both of $A\in\mathcal M$, $B\in\mathcal M$ hold. Let us say that $A\notin\mathcal M$ (the proof is entirely symmetric if we assume $B\notin\mathcal M$ instead). Then by Lemma~\ref{lem:annexion}, there is a $C\in\mathcal M$ that is disjoint from $A$, which in turn implies that $C\subseteq X\setminus A$. Since $X\setminus A\subseteq X\setminus(A\cap B)$ and $\mathcal M$ is closed upwards by Lemma~\ref{lem:annexion}, we conclude that $X\setminus(A\cap B)\in\mathcal M$, and we are done. 
\end{proof}

\begin{remark}
Ultrafilters are usually defined as maximal filters, where a filter is defined to be a family of subsets of $X$ that is nonempty, closed upwards and closed under intersections, and does not contain $\varnothing$. So the only difference between a filter and a nice family as defined here is the requirement that the family be closed upwards. However, by Lemma~\ref{lem:annexion}, this distinction vanishes once we focus on maximal objects (that is, being a maximal filter is equivalent to being a maximal nice family), and so by Lemma~\ref{lem:maximalnice}, we see that our (somewhat nonstandard\footnote{The author first came in contact with this nonstandard definition in a set of notes by Andreas Blass.}) definition of an ultrafilter is equivalent to the most common one.
\end{remark}

\begin{corollary}
If $X$ is an arbitrary infinite set, then there exists a nonprincipal ultrafilter on $X$.
\end{corollary}

\begin{proof}
Given such an $X$, notice that the family of all cofinite subsets of $X$,
\begin{equation*}
\fr=\{A\subseteq X\mid X\setminus A\text{ is finite}\}
\end{equation*}
(usually known as the {\em Fr\'echet filter on $X$}), is a nice family. Hence, by Zorn's lemma, there exists an ultrafilter $u$ with $\fr\subseteq u$. Clearly $u$ is nonprincipal (otherwise there would be a finite $F$ with $F\in u$ by Theorem~\ref{thm:principaluf}, and since $X\setminus F\in\fr\subseteq u$, this would mean that $\varnothing=F\cap(X\setminus F)\in u$, a contradiction).
\end{proof}

\begin{remark}
Even though in this article we used Zorn's lemma to prove the existence of nonprincipal ultrafilters, we do not actually need the full strength of the axiom of choice for this. In fact, the statement that for every infinite $X$ there exists a nonprincipal ultrafilter on $X$ is equivalent over the $\zf$ axioms to the Boolean prime ideal theorem (which is in turn equivalent to Tychonoff's theorem for Hausdorff topological spaces, and strictly weaker than the full axiom of choice). See~\cite{howard-rubin} for an extensive compendium of statements weaker than the axiom of choice and their multiple relations of implication and equivalence.
\end{remark}

\subsection{Application: Ramsey's theorem.}

A paradigmatic example of a Ramsey-theoretic result is the statement that in every party with at least six attendees, one can find either three of them that do not know each other, or three of them that mutually know each other. In mathematical terms, we can model the situation by means of a complete graph with as many vertices as attendees at the party, together with a coloring of the edges of this graph---say, we declare an edge to be \tr\ if the two extremes of that edge are vertices representing two people who know each other, and we make the edge \tb\ otherwise. Then, the result just mentioned states that, regardless of what the coloring of the edges is, we will always be able to find a triangle that is either all \tr\, or all \tb: we say that we can always find a {\em monochromatic triangle}.

This is a recurring theme within Ramsey theory: one colors a structure (mathematically, colorings are represented by functions whose codomain is some finite set; in this article, we will use the set $\{\r,\b\}$),\footnote{Throughout this article we consider colorings with only two colors in order to simplify matters, since for the problems that we study here this does not make a difference, and the same arguments work with any finite number of colors. However, the reader should be warned that this is not the case in certain contexts, as there are examples of Ramsey-type properties that are true for 2-colorings but fail for arbitrary finite colorings~\cite{csikvari-gyarmati-sarkozy}. As a particularly striking recent example, for every coloring of $\mathbb N$ with two colors, there are infinitely many monochromatic solutions to the Diophantine equation $x+y=z^2$, whereas there exists a coloring of $\mathbb N$ with three colors such that no nontrivial solution to the same equation can be monochromatic~\cite{green-lindqvist}.} and if the structure is large enough, then one can usually find some sufficiently rich monochromatic substructures.

Let us prove the weaker version of the aforementioned theorem where we require the party to contain infinitely many attendees.

\begin{theorem}
Let $G=(V,E)$ be an infinite complete graph. For every coloring $c:E\longrightarrow\{\r,\b\}$, there are three vertices $x,y,z$ such that the subgraph of $G$ induced by $\{x,y,z\}$ is monochromatic.
\end{theorem}

\begin{proof}
Let $u$ be a nonprincipal ultrafilter on $V$. Upon fixing a vertex $v\in V$, notice that we have a partition of $V$
\begin{equation*}
V=\{v\}\cup\{w\in V\mid c(vw)=\r\}\cup\{w\in V\mid c(vw)=\b\}.
\end{equation*}
Exactly one of these pieces belongs to $u$ by Proposition~\ref{prop:finiteunion}. Since $u$ is nonprincipal, $\{v\}\notin u$ so either $\{w\in V\mid c(vw)=\r\}\in u$, in which case we will say that $v$ is $u$-\tr, or $\{w\in V\mid c(vw)=\b\}\in u$, in which case we will say that $v$ is $u$-\tb.

The procedure described in the previous paragraph can, of course, be carried out for every vertex $v\in V$; so each $v\in V$ is either $u$-\tr\ or $u$-\tb. This induces yet another partition of $V$, given by
\begin{equation*}
V=\{v\in V\mid v\text{ is }u\text{-\tr}\}\cup\{v\in V\mid v\text{ is }u\text{-\tb}\},
\end{equation*}
and exactly one piece of this partition belongs to $u$, again by Proposition~\ref{prop:finiteunion}. Let us assume, without loss of generality, that $A=\{v\in V\mid v\text{ is }u\text{-\tr}\}\in u$, and construct a monochromatic triangle in color \tr---otherwise, we would be able to construct a monochromatic triangle in color \tb, in an entirely symmetric fashion. For each $v\in A$, the fact that $v$ is $u$-\tr\ means that $u$ contains the set $A_v=\{w\in V\mid c(vw)=\r\}$. We begin by taking $v_1\in A$ and $v_2\in A\cap A_{v_1}$. We can do this because, since $A\in u$, we have that $A$ is nonempty; now $v_1\in A$ and so $A_{v_1}\in u$, therefore we have $A\cap A_{v_1}\in u$ (because ultrafilters are closed under finite intersections) and in particular $A\cap A_{v_1}$ is nonempty. Since $v_2\in A_{v_1}$, we have $c(v_1 v_2)=\r$. Now, since $v_2\in A$, we have $A_{v_2}\in u$ and so we can pick a $v_3\in A_{v_1}\cap A_{v_2}$, because $A_{v_1},A_{v_2}\in u$ and consequently $A_{v_1}\cap A_{v_2}\in u$; in particular, $A_{v_1}\cap A_{v_2}$ is nonempty. Since $v_3\in A_{v_1}$, we have $c(v_1 v_3)=\r$; since $v_3\in A_{v_2}$, we have $c(v_2 v_3)=\r$. We already argued that $c(v_1 v_2)=\r$ as well; it follows that the triangle formed by the vertices $v_1,v_2,v_3$ is monochromatic in color \tr, and we are done.
\end{proof}

The reader might rightly complain that we did not prove the theorem that we announced at the beginning---the one where the graph $G$ only has 6 vertices---but rather a significantly weaker version of it. In our defense, we will point out that, by proving this weaker version of the theorem the way we did, we can now readily improve the theorem, not by reducing the number of vertices needed from the graph $G$, but by increasing the size of the monochromatic structure. To see this, notice that, at the moment when we picked $v_3$ in the proof above, we could have picked it to be an element of $A\cap A_{v_1}\cap A_{v_2}$ (since ultrafilters are closed under finite intersections, this set is an element of $u$, and hence nonempty). Then, the fact that $v_3\in A_{v_1}\cap A_{v_2}$ still ensures that the triangle formed by the vertices $v_1,v_2,v_3$ is monochromatic in color \tr, whereas the fact that $v_3\in A$ implies that, furthermore, we also have $A_{v_3}\in u$. Hence, we can now also pick yet another vertex $v_4\in A\cap A_{v_1}\cap A_{v_2}\cap A_{v_3}$ (this set is nonempty by virtue of its belonging to $u$), and see that, for each $i\in\{1,2,3\}$, $v_4\in A_{v_i}$ and therefore $c(v_i v_4)=\r$. Thus we have obtained vertices $v_1,v_2,v_3,v_4$ such that all the edges between these four vertices are \tr; we say that the (complete) subgraph of $G$ induced by $\{v_1,v_2,v_3,v_4\}$ is monochromatic in color \tr; note that, furthermore, we have $v_4\in A$.

We can continue the process described in the previous paragraph for as long as we wish: reasoning inductively, if we already have $n$ vertices $v_1,\ldots,v_n$ such that $A_{v_1},\ldots,A_{v_n}\in u$ and such that the (complete) subgraph of $G$ induced by $\{v_1,\ldots,v_n\}$ is monochromatic in color \tr, use the fact that ultrafilters are closed under finite intersections to pick a further vertex $v_{n+1}\in A\cap A_{v_1}\cap\cdots\cap A_{v_n}$ (since the latter set belongs to $u$ and hence is nonempty). For each $i\in\{1,\ldots,n\}$, the fact that $v_{n+1}\in A_{v_i}$ means that $c(v_i v_{n+1})=\r$; this, together with the fact that all the edges between $\{v_1,\ldots,v_n\}$ are \tr\ implies that the complete subgraph of $G$ induced by the vertices $\{v_1,\ldots,v_n,v_{n+1}\}$ is also monochromatic in color \tr. On the other hand, the fact that $v_{n+1}\in A$ implies that $A_{v_{n+1}}\in u$ and so the induction can continue for one more step. Hence, we can obtain arbitrarily large (complete) subgraphs of $G$ that are monochromatic; in fact, by following through the induction for all the natural numbers, one actually obtains an infinite complete subgraph of $G$ (the one induced by the set of vertices $\{v_n\mid n\in\mathbb N\}$) that is monochromatic. We record this fact in the theorem below.

\begin{theorem}[Ramsey~\cite{ramsey}]
Let $G=(V,E)$ be an infinite complete graph. For every coloring $c:E\longrightarrow\{\r,\b\}$, there exists an infinite subset $X\subseteq V$ such that the subgraph of $G$ induced by $X$ is monochromatic.
\end{theorem}

To finish the section, we wish to point out a couple of important concepts in Ramsey theory, namely {\it partition regularity} and {\it weak partition regularity}. An upwards closed family $\mathcal F$ of subsets of $\mathbb N$ is called {\it weakly partition regular} if for every finite partition of $\mathbb N$, one of the cells belong to $\mathcal F$, and it is called {\it partition regular} if for every partition of a member of $\mathcal F$, one of the cells belongs to $\mathcal F$. In an abstract sense, it can be said that Ramsey theory consists of the study of weak partition regular and partition regular properties. The reader is encouraged to work out the translation between the Ramsey-theoretic statements proved in this article and the weak partition regularity of certain families of subsets of $\mathbb N$ (in fact, for all the statements addressed in this article, the corresponding family of sets will be partition regular). Regarding Ramsey-theoretic statements in this light makes the connection with ultrafilters very explicit, since a family $\mathcal F$ is weakly partition regular if and only if there is an ultrafilter $u$ with $u\subseteq\mathcal F$~\cite[Theorem 5.7]{hindman-strauss}, and it is partition regular if and only if $\mathcal F$ is a union of ultrafilters~\cite[Theorem 3.11]{hindman-strauss}.

\section{Additive structure.}

We will now introduce the algebraic-topological structure of the collection of all ultrafilters on a given set, and then illustrate how the theory about this structure allows us to prove some Ramsey-theoretic results in additive combinatorics. Since we are mostly interested in illustrating the applications of the theory, rather than in presenting the theory itself, in the first subsection we will for the most part simply mention the relevant results without proofs.\footnote{The reader interested in delving deep into the development of this theory can consult~\cite{hindman-strauss}, or, for a rather compact, but complete, exposition of all the tools used in this article, see~\cite[Section 2.1, pp. 27--30]{stevo-ramsey}.} In the following subsection we will, on the other hand, fully explain the applications of these results with all the details.

\subsection{Algebra in the \v{C}ech--Stone compactification.}

Given a set $X$, we denote the set of all ultrafilters on $X$ by $\beta X$. The set $\beta X$ can be topologized by declaring, for each $A\subseteq X$, the set
\begin{equation*}
\overline{A}=\{u\in\beta X\mid A\in u\}
\end{equation*}
to be a basic open set. With this topology (known as the \emph{Stone topology}), the space $\beta X$ is compact Hausdorff, known as the \textbf{\v{C}ech--Stone compactification} of $X$, and it contains a dense copy of $X$ via the embedding $x\longmapsto u_x$ that maps each point $x\in X$ to the corresponding principal ultrafilter $u_x$. Once we identify $X$ with its copy within $\beta X$ (such a copy is a discrete subspace of $\beta X$, so we can think of $X$ as a discrete space that is embedded into $\beta X$), we can see that $\overline{A}$ is a clopen set which is really the closure of $A\subseteq X$ within $\beta X$. The \textbf{\v{C}ech--Stone remainder} of $X$ is the closed (and hence compact) subspace of $\beta X$ consisting of all nonprincipal ultrafilters on $X$; this subspace is usually denoted $X^*=\beta X\setminus X$. Given two sets $X,Y$, every function $f:X\longrightarrow Y$ lifts to a (unique) continuous extension (still denoted by $f$, and known as the \textit{\v{C}ech--Stone extension} of $f$) $f:\beta X\longrightarrow\beta Y$, which is given by $f(u)=\{A\subseteq Y\mid f^{-1}[A]\in u\}$. The ultrafilter $f(u)$ is sometimes known as the \textit{Rudin--Keisler image} of the ultrafilter $u$; this ultrafilter can also be described as the ultrafilter on $Y$ generated by the family $\{f[A]\mid A\in u\}$.


The previous paragraph describes the topological structure of $\beta X$, and we will now proceed to describe the corresponding algebraic structure. Suppose that, rather than a bare set $X$ without any additional structure, we have a semigroup operation $*$ on $X$. We can now define a semigroup operation on $\beta X$ by means of the formula
\begin{equation*}
u*v=\{A\subseteq X\mid \{x\in X\mid \{y\in X\mid x*y\in A\}\in v\}\in u\};
\end{equation*}
this operation turns $\beta X$, equipped also with the Stone topology, into a compact right-topological semigroup. This means that $\beta X$ is a semigroup when equipped with the operation $*$ and, for each fixed $v\in\beta X$, the right translation mapping $u\longmapsto u*v$ is a continuous function from $\beta X$ to $\beta X$ (although left translations $v\longmapsto u*v$ need not be continuous, and the semigroup operation $*$ need not be jointly continuous). Furthermore, if the semigroup $X$ is sufficiently well-behaved,\footnote{Technically, the condition needed on the semigroup $X$ is that for every finite $F\subseteq X$ and for every infinite $A\subseteq X$ there are finitely many $a_1,\ldots,a_n\in A$ such that the set $\{x\in X\mid a_i*x\in F\text{ for all }i\}$ is finite.\label{foot:semigroup}} then the set $X^*=\beta X\setminus X$ of nonprincipal ultrafilters on $X$ is a closed subsemigroup of $\beta X$, and hence $X^*$ is a compact right-topological semigroup in its own right. Whenever we have semigroup operations on both sets $X$ and $Y$, and a function $f:X\longrightarrow Y$ is a semigroup homomorphism (meaning $f(x*y)=f(x)*f(y)$ for all $x,y\in X$), the corresponding \v{C}ech--Stone extension $f:\beta X\longrightarrow\beta Y$ will be a continuous semigroup homomorphism.

The definition of the semigroup operation on $\beta X$ might not look very natural at first sight, but it in fact arises naturally in various different ways. For example, remember that when we first introduced ultrafilters we mentioned that they can be thought of as (finitely additive) measures on the set $X$. Say that we have two ultrafilters $u,v\in\beta X$ and we use the letters $\mu,\nu$ to denote the respective measures (thus $\mu(A)=1$ if and only if $A\in u$, and $\mu(A)=0$ otherwise; similarly $\nu(A)=1$ if and only if $A\in v$ and $\nu(A)=0$ otherwise). Then the measure $\mu\diamond\nu$ that corresponds to the ultrafilter $u*v$ is simply given by
\begin{equation*}
(\mu\diamond\nu)(A)=\int_{x\in X}\int_{y\in X}\chi_A(x*y)d\nu(y) d\mu(x)
\end{equation*}
(where $\chi_A$ is the characteristic function of $A$), for every $A\subseteq X$. Thus, we can see that, from the appropriate viewpoint, the semigroup operation $u*v$ behaves like a convolution.

To finish this section, we will prove a crucial lemma that is of central importance for the algebra in the \v{C}ech--Stone compactification. We first introduce the relevant terminology.

\begin{definition}\hfill
\begin{enumerate}
\item Let $(S,*)$ be a semigroup. The element $x\in S$ is said to be {\bf idempotent} if it satisfies $x*x=x$.
\item A triple $(S,*,\tau)$ is a {\bf right-topological semigroup} if $(S,*)$ is a semigroup and $(S,\tau)$ is a Hausdorff topological space such that, for every $s\in S$, the right-translation mapping $S\longrightarrow S$ given by $t\longmapsto t*s$ is continuous.
\end{enumerate}
\end{definition}

Hence, one of the crucial facts that we outlined in this section is that, whenever $X$ is a set with a semigroup operation $*$, the triple $(\beta X,*,\tau)$ is a compact right-topological semigroup (where $*$ denotes the extension of the operation $*:X\times X\longrightarrow X$ to all of $\beta X$ and $\tau$ denotes the Stone topology on $\beta X$).

\begin{lemma}[Ellis--Numakura~\cite{ellis-lemma,numakura-lemma}]\label{ellis-numakura}
If $(S,*,\tau)$ is a compact right-topological semigroup, then there exists an idempotent element $s\in S$.
\end{lemma}

\begin{proof}
To begin with, we remark that, since $S$ is compact and Hausdorff, every subset $K\subseteq S$ is closed if and only if it is compact, and so we will utilize the words ``closed'' and ``compact'' interchangeably throughout the proof.

We start by noting that the collection of all nonempty $T\subseteq S$ that are subsemigroups of $S$ (that is, closed under the operation $*$) and compact, ordered by reverse inclusion, satisfies the hypotheses of Zorn's lemma (for any $\subseteq$-linearly ordered collection of compact subsemigroups of $S$, its intersection is also a compact subsemigroup---furthermore, nonempty, by compactness of $S$) and therefore there exists a minimal closed subsemigroup $T\subseteq S$. Let $s\in T$ be arbitrary. Then $T*s\subseteq T$, and it is readily checked that $T*s$ is a subsemigroup of $S$ which is, furthermore, compact (as it is the continuous image of the compact set $T$ under the continuous mapping $t\longmapsto t*s$). Hence, by minimality of $T$, we must have $T*s=T$; in particular, $s\in T*s$. This shows that the set $\{t\in T\mid t*s=s\}$ is nonempty. Furthermore, the set $\{t\in T\mid t*s=s\}\subseteq T$ is compact (as it is the preimage of the closed set $\{s\}$ under the continuous mapping $t\longmapsto t*s$), and it can be readily checked that it is a subsemigroup of $S$; so by minimality of $T$, we must have $T=\{t\in T\mid t*s=s\}$ and in particular, since $s\in T$, we have $s*s=s$, which shows that $s$ is idempotent and finishes the proof.
\end{proof}

\begin{corollary}\label{cor:existidemp}
Let $(X,*)$ be a (sufficiently well-behaved, see footnote~\ref{foot:semigroup}) semigroup. Then there are idempotent elements (that is, ultrafilters $u$ such that $u*u=u$) in $\beta X$. Furthermore, whenever $Z\subseteq\beta X$ is a closed subsemigroup, one can find an idempotent ultrafilter $u\in Z$. In particular, there are nonprincipal idempotent ultrafilters on $X$.
\end{corollary}

\begin{proof}
This is immediate from the fact that $\beta X$ is a compact right-topological semigroup (in particular, any closed subset of it is also compact) and from Lemma~\ref{ellis-numakura}. The last statement follows from the fact that the set $X^*=\beta X\setminus X$ of nonprincipal ultrafilters on $X$ is a closed subsemigroup of $\beta X$.
\end{proof}

Suppose that we have a semigroup $X$ and an idempotent ultrafilter $u\in\beta X$. Using the definition of the semigroup operation $*$ on $\beta X$, we see that, for every $A\in u$, it is the case that (since $u*u=u$ and thus $A\in u*u$)
\begin{equation*}
\{x\in X\mid \{y\in X\mid x*y\in A\}\in u\}\in u.
\end{equation*}
This observation will be extremely important in the applications that we illustrate in the following subsection.

\subsection{Application: Schur, Folkman--Rado--Sanders, and Hindman's theorems.}

Leaving graph theory behind, we can also encounter a large number of Ramsey-theoretic results in the field of additive combinatorics; as the name suggests, we are talking here about results where the monochromatic structures obtained are defined, in one way or another, in terms of the addition operation. Possibly the oldest result in this vein is Schur's theorem~\cite{schur}, which establishes that, whenever we color the elements of $\mathbb N$ with finitely many colors (that is, whenever we partition the set $\mathbb N$ into finitely many cells), it is always possible to find two distinct elements $x,y$ such that the set $\{x,y,x+y\}$ is monochromatic. With the theory presented so far, we are in a position to state and prove this theorem right away.

\begin{theorem}[Schur~\cite{schur}]\label{thm:schur}
Suppose that we have a coloring $c:\mathbb N\longrightarrow\{\r,\b\}$. Then there exists a monochromatic set of the form $\{x,y,x+y\}$.
\end{theorem}

\begin{proof}
Use the Ellis--Numakura lemma (or rather, Corollary~\ref{cor:existidemp}) to obtain a nonprincipal idempotent ultrafilter $u$ on $\mathbb N$. Since $\mathbb N=\{n\mid n\text{ is \tr}\}\cup\{n\mid n\text{ is \tb}\}$, we have that either $\{n\mid n\text{ is \tr}\}\in u$ or $\{n\mid n\text{ is \tb}\}\in u$; suppose without loss of generality (the proof being entirely symmetric otherwise) that $A=\{n\mid n\text{ is \tb}\}\in u$. Since $u$ is an idempotent, we know that 
\begin{equation*}
A^\star=A\cap\{n\in\mathbb N\mid \{m\in\mathbb N\mid n+m\in A\}\in u\}\in u.
\end{equation*}
Take an $x\in A^\star$. Then $x$ is \tb, and moreover, the set
\begin{equation*}
A_x=\{m\in\mathbb N\mid x+m\in A\}\in u,
\end{equation*}
so we can also take a $y\in A_x\cap A$. Thus, we will have that (since $y\in A$) $y$ is also \tb, and, moreover, since $y\in A_x$, we have that $x+y\in A$ and thus $x+y$ is also \tb. Hence, we have shown that the set $\{x,y,x+y\}$ consists entirely of \tb\ elements.
\end{proof}

Notice how, in the previous proof, the fact that the set $A$ is precisely the set of all elements that were colored \tb\ is not important for the proof itself. If we erase all references to colors from the previous proof, and instead just work with an arbitrary set $A$ satisfying that $A\in u$, everything works out just as nicely, and we obtain a useful observation.

\begin{remark}\label{rmk:schur}
Our proof of Theorem~\ref{thm:schur} actually shows the following: whenever $u$ is a nonprincipal idempotent ultrafilter on $\mathbb N$, for any $A\in u$ we can find two distinct elements $x,y$ such that $\{x,y,x+y\}\subseteq A$.
\end{remark}

Using Remark~\ref{rmk:schur}, we can improve Theorem~\ref{thm:schur} to obtain even more: Suppose that $\mathbb N$ is colored with colors \tb\ and \tr, take a nonprincipal idempotent ultrafilter $u$ over $\mathbb N$, and assume without loss of generality that $A=\{n\in\mathbb N\mid n\text{ is }\b\}\in u$. Since $u$ is idempotent, we have $\{n\in A\mid \{m\in A\mid n+m\in A\}\in u\}\in u$ and so we can choose an $x$ belonging to that set; this means both that $x$ is \tb\ and that the set $B=\{m\in\mathbb N\mid x+m\text{ is }\b\}\in u$. By Remark~\ref{rmk:schur}, we can now find two elements $y,z$ such that $y,z,y+z\in A\cap B$ (because, since ultrafilters are closed under finite intersections, we have $A\cap B\in u$). This means both that $y,z,y+z$ are \tb\ (since they are elements of $A$), and also (by the definition of $B$) that $x+y,x+z,x+y+z$ are \tb. We have thus concluded that the entire set $\{x,y,z,x+y,x+z,y+z,x+y+z\}$ consists of \tb\ elements (or, had $u$ contained the set of \tr\ elements instead, we would have obtained such a set consisting of \tr\ elements). This means that for every coloring of $\mathbb N$ with two colors, it is possible to obtain three distinct elements $x,y,z$ such that the corresponding {\it set of all finite sums} $\{x,y,z,x+y,x+z,y+z,x+y+z\}$ is monochromatic. As we may expect by now, this theorem can be significantly generalized, and we will now show how.

\begin{definition}
Given a set of numbers $X\subseteq\mathbb N$, the set of {\bf finite sums} of $X$ is defined to be
\begin{eqnarray*}
\fs(X) & = & \left\{x_1+\cdots+x_n\mid n\in\mathbb N,\ x_1,\ldots,x_n\text{ are distinct and }x_1,\ldots,x_n\in X\right\} \\
 & = & \left\{\sum_{x\in F}x\mid F\subseteq X\text{ is finite nonempty}\right\}.
\end{eqnarray*}
\end{definition}

\begin{theorem}[Folkman and Rado~
(see~\cite{graham-rothschild-spencer}, Thms. 3.5 and 3.4)
, Sanders~\cite{sanders}]\label{folkman-rado-sanders}
Given a coloring of $\mathbb N$ with colors \tr\ and \tb, and given an $n\in\mathbb N$, it is possible to find elements $x_1,\ldots,x_n\in\mathbb N$ such that the set $\fs(x_1,\ldots,x_n)$ is monochromatic.
\end{theorem}

\begin{proof}
Let $u$ be a nonprincipal idempotent ultrafilter on $\mathbb N$. We proceed to show, by induction on $n\in\mathbb N$, that if $A\in u$ then there are $n$ distinct elements $x_1,\ldots,x_n$ such that $\fs(x_1,\ldots,x_n)\subseteq A$. The case $n=1$ is obvious, so suppose the statement holds for $n$, and let $A\in u$. Since $u$ is an idempotent element, we have $A^\star=A\cap\{x\in\mathbb N\mid \{y\in\mathbb N\mid x+y\in A\}\in u\}\in u$. In particular $A^\star$ is nonempty and so we can choose an $x_1\in A^\star$. This means that $x_1\in A$ and furthermore $B=A\cap\{y\in\mathbb N\mid x_1+y\in A\}\in u$; so by induction hypothesis we can find distinct elements $x_2,\ldots,x_{n+1}$ such that $\fs(x_2,\ldots,x_{n+1})\subseteq B$. Note that, if $z\in\fs(x_1,\ldots,x_{n+1})$, then either $z\in\fs(x_2,\ldots,x_{n+1})\subseteq B\subseteq A$, or $z=x_1+y$ for some $y\in\fs(x_2,\ldots,x_{n+1})\subseteq B$, which means, by definition of $B$, that $z=x_1+y\in A$. In either case we get $z\in A$, which implies that $\fs(x_1,\ldots,x_{n+1})\subseteq A$, finishing the induction.

Now suppose that $\mathbb N$ has been colored with colors \tb\ and \tr, and let $n\in\mathbb N$. Exactly one of the two sets $\{x\in\mathbb N\mid x\text{ is }\b\}$ and $\{x\in\mathbb N\mid x\text{ is }\r\}$ belongs to $u$; let us denote whichever set this is by $A$. Then, by the statement proved in the previous paragraph, there are elements $x_1,\ldots,x_n$ such that $\fs(x_1,\ldots,x_n)\subseteq A$. By choice of $A$, this means precisely that $\fs(x_1,\ldots,x_n)$ is monochromatic.
\end{proof}

A more careful look at the proof of Theorem~\ref{folkman-rado-sanders} reveals that it should, in fact, be possible to use the same ideas to prove an even stronger result, in which one finds {\it infinitely many} elements whose set of finite sums is monochromatic. This result is known as Hindman's theorem, and its history is quite interesting. What is now known as Hindman's theorem was, for a relatively long time, a conjecture of Graham and Rothschild (asked as a question in~\cite[p. 291]{graham-rothschild}, and attributed to those authors as a conjecture by Erd\H{o}s in~\cite[p. 122]{erdos}); it was eventually proved by Hindman using involved combinatorial arguments, and subsequently Baumgartner~\cite{baumgartner-short-proof-of-hindman} found a somewhat simpler, but still purely combinatorial proof. However, Galvin, even before Hindman's proof (back when Hindman's theorem was still just a conjecture), had realized that the existence of an ultrafilter with a certain combinatorial property would readily imply the theorem, but he was unable to prove that these special ultrafilters exist. (In retrospect, we know that the ultrafilters that Galvin considered are precisely the idempotent ultrafilters.) Hindman~\cite{hindman-ultrafiltersch} proved that these special ultrafilters exist if one assumes both the continuum hypothesis and the conjecture that would eventually become Hindman's theorem, but the question whether these ultrafilters {\it actually exist} (i.e., whether their existence can be established without additional assumptions) remained open---even after Hindman's theorem had become an actual theorem. Later, one day Galvin asked Glazer whether these ultrafilters that he had considered existed, and, when Glazer answered affirmatively almost instantly and without hesitation, Galvin's reaction was one of disbelief, thinking that Glazer must have misunderstood the question, which could not have been that easy to answer. It turns out that there was, in fact, no misunderstanding, and the question had indeed a simple answer: at the time, very few people working on the algebra in the \v{C}ech--Stone compactification of $\mathbb N$ were aware that such compactification can be seen as a space of ultrafilters, but Glazer did know this, and so he knew that $\beta N$ has idempotent elements, and it was not hard for him to see that the combinatorial property that Galvin was seeking was in fact equivalent to idempotence of the relevant ultrafilter. Neil Hindman has told this story in various places,\footnote{The story, as presented here, seems to agree with Galvin's own recollections, as verified by the author through personal communication.} such as~\cite[pp. 120--121]{hindman-number-theory}, \cite[pp. 835--836]{hindman-steprans-strauss}, or~\cite[pp. 122--123]{hindman-strauss}. We will now finally proceed to prove this outstanding result (the proof below is the Galvin--Glazer one\footnote{There is also a proof of Hindman's theorem, based on tools from topological dynamics, due to Furstenberg and Weiss~\cite{furstenberg-weiss}.}).

\begin{theorem}[Hindman~\cite{hindmanhisthm}]
Given any coloring $c:\mathbb N\longrightarrow\{\r,\b\}$, there exists an infinite set $X=\{x_1,\ldots,x_n,\ldots\}\subseteq\mathbb N$ such that the set $\fs(X)$ is monochromatic.
\end{theorem}

\begin{proof}[Proof (sketch)]
Let us analyze further the proof of the Folkman--Rado--Sanders theorem, Theorem~\ref{folkman-rado-sanders} above. There, we begin by taking a nonprincipal idempotent ultrafilter $u$ over $\mathbb N$, fix an arbitrary $A\in u$, and proceed to prove by induction on $n\in\mathbb N$ that there are elements $x_1,\ldots,x_n$ such that $\fs(x_1,\ldots,x_n)\in A$.

Now, when running the induction, the inductive step starts in exactly the same way, regardless of the value of $n$: choose any $x_1\in A^\star=A\cap\{x\in\mathbb N\mid \{y\in\mathbb N\mid x+y\in A\}\in u\}$. From there, the induction hypothesis allows us to find $x_2,\ldots,x_{n+1}$ such that $\fs(x_2,\ldots,x_{n+1})\subseteq A\cap\{y\in\mathbb N\mid x_1+y\in A\}$; this is what ensures that $\fs(x_1,x_2,\ldots,x_{n+1})\subseteq A$. Notice that, for whatever value of $x_1$ that we fix at the beginning (as long as it belongs to $A^\star$), the remainder of the proof works equally well for every $n$. So we can see that the following statement holds:

\begin{center}
For any $n\geq 1$, there are $x_1,\ldots,x_n$ such that $\fs(x_1,\ldots,x_n)\subseteq A$,\\ {\it and the value of $x_1$ does not depend on $n$}.
\end{center}

In other words, it is possible to find a single number $x_1\in\mathbb N$ such that for every $n$, $x_1$ is the first element of a set $X$ with $n$ elements satisfying $\fs(X)\subseteq A$. But even more is true. In fact, going back once again to the proof of Theorem~\ref{folkman-rado-sanders} above, we see that, once $x_1\in A^\star$ has been fixed, we obtain the remaining elements $\{x_2,\ldots,x_n\}$ satisfying $\fs(x_1,x_2,\ldots,x_n)$ by choosing them in such a way that $\fs(x_2,\ldots,x_n)\subseteq B=A\cap\{y\in\mathbb N\mid x_1+y\in A\}$ (which is the induction hypothesis in action). The attentive reader will notice that we can now apply the same reasoning as above to the set $B$, in order to obtain $x_2$ (that is, pick any $x_2\in B^\star$). So we see that with this reasoning we can conclude the following:

\begin{center}
For any $n\geq 3$, there are $x_1,\ldots,x_n$ such that $\fs(x_1,\ldots,x_n)\subseteq A$,\\ {\it and the values of $x_1$ and $x_2$ do not depend on $n$}.
\end{center}

The reader should realize now that this reasoning can be carried out over and over. Thus, the following statement (which the interested reader is encouraged to prove formally by herself, using induction) should ring true, for every $k\in\mathbb N$:

\begin{center}
For every $n>k$, there are elements $x_1,\ldots,x_n$ such that $\fs(x_1,\ldots,x_n)\subseteq A$,\\ {\it and the values of $x_1,\ldots,x_k$ do not depend on $n$}.
\end{center}

With the previous claim under our belt, we proceed as follows: By recursively fixing the values of the $x_k$ encountered in the previous reasoning, we obtain an infinite sequence of elements $x_1,\ldots,x_n,\ldots,$ satisfying that, for every $n\in\mathbb N$, it is the case that $\fs(x_1,\ldots,x_n)\subseteq A$. Now, it is not hard to see that this in fact implies that, if $X=\{x_1,\ldots,x_n,\ldots\}$, then $\fs(X)\subseteq A$: for if $y\in\fs(X)$, then there is a sufficiently large $n\in\mathbb N$ such that $y\in\fs(x_1,\ldots,x_n)$, and so $y\in A$. The conclusion is that, given any nonprincipal idempotent ultrafilter $u$ and an element $A\in u$, there is an infinite set $X$ such that $\fs(X)\subseteq A$.

Now suppose that we have a coloring $c:\mathbb N\longrightarrow\{\r,\b\}$. Let $u$ be a nonprincipal idempotent ultrafilter on $\mathbb N$, and let $A\in u$ be the set containing either all \tr\ elements, or all \tb\ elements ($u$ contains one and only one of these two possibilities). Then, by the statement whose proof we just sketched, one can find an infinite set $X\subseteq\mathbb N$ such that $\fs(X)\subseteq A$; in other words, $\fs(X)$ is monochromatic (in one color or the other, depending on what color the set $A$ represents).
\end{proof}

\section{Another application to additive combinatorics.}

In order to present our third and last application of the theory of ultrafilters to obtain Ramsey-theoretic results, we will need to introduce some more theory regarding the algebraic-topological structure of the semigroup of ultrafilters. We will proceed in largely the same way as in the previous section: first we will present the additional ultrafilter-theoretic results without proofs,\footnote{As noted before, an interested reader can consult~\cite[Section 2.1, pp. 27--30]{stevo-ramsey} for a short but complete account of all the results mentioned here, or look into~\cite{hindman-strauss} for a more extensive treatment.} and then we will present full details of how these results can be applied to provide an ultrafilter proof of van der Waerden's theorem.

\subsection{Minimal idempotents and ideals.}

We proceed to introduce the notions of minimal idempotents and ideals.

\begin{definition}
Let $(S,*)$ be a semigroup.
\begin{enumerate}
\item We define $E(S)=\{x\in S\mid x\text{ is an idempotent}\}$,
\item We define a partial order relation $\leq$ on $E(S)$ by stipulating that $x\leq y$ if and only if $x*y=y*x=x$.
\item A subset $I\subseteq S$ is an {\em ideal} of $S$ if for every $x\in S$, $y\in I$, we have $x*y\in I$ and $y*x\in I$.
\end{enumerate}
\end{definition}

Note that, if we have two semigroups $(S,*)$ and $(T,*)$, a semigroup homomorphism $f:S\longrightarrow T$, and an idempotent element $x\in E(S)$, then (as $f(x)*f(x)=f(x*x)=f(x)$) we have that $f(x)$ is an idempotent element of $T$. Furthermore, if $x,y\in E(X)$ and $x\leq y$, then $f(x)\leq f(y)$ (since $f(x)*f(y)=f(x*y)=f(y*x)=f(y)*f(x)=f(y*x)=f(x)$). Thus, semigroup homomorphisms preserve idempotent elements, as well as the partial order relation just defined between them.

Now suppose that we have a compact right-topological semigroup $(S,*,\tau)$. A crucial result in this context is that minimal idempotents (that is, idempotent elements of $S$ that are minimal with respect to the partial order $\leq$) exist. This can be proved using Zorn's lemma together with Ellis's lemma, along with a result ensuring that an idempotent element of $S$ is minimal if and only if it belongs to some minimal closed left ideal.\footnote{A closed left ideal of $S$ is a closed subset $I\subseteq S$ satisfying $x*y\in I$ whenever $x\in S$ and $y\in I$.} Furthermore, whenever we have an arbitrary idempotent element $x\in E(S)$ and an arbitrary closed ideal $J\subseteq S$, we can always find a minimal idempotent $y\in J$ with $y\leq x$. This last statement will be of central importance for the application that we will illustrate in the next subsection.

The results mentioned in the previous paragraph will be used in the context where the compact right-topological semigroup that we deal with is the \v{C}ech--Stone compactification $\beta X$ of some semigroup $X$. It is possible to verify (though we will not do so here) that, if $X$ is a semigroup and $Y\subseteq X$ is a subsemigroup, then $\overline{Y}=\{u\in\beta X\mid Y\in u\}$ is a closed subsemigroup of $\beta X$. In fact, if $Y\subseteq X$ is infinite, then $\{u\in\beta X\mid u\text{ is nonprincipal and }Y\in u\}$ is also a closed subsemigroup of $\beta X$ and therefore it contains idempotent elements by Lemma~\ref{ellis-numakura}, so we can always ensure the existence of nonprincipal idempotent ultrafilters belonging to $\overline{Y}$. Similarly, if $I\subseteq X$ is an ideal of $X$, then $\overline{I}=\{u\in\beta X\mid I\in u\}$ is a closed ideal of $\beta X$, and so is the set $\{u\in\beta X\mid u\text{ is nonprincipal and }I\in u\}$ whenever $I$ is infinite. In particular, if $u\in\beta X$ is an idempotent element, then one can find a nonprincipal minimal idempotent $v\in\overline{I}$ with $v\leq u$.

\subsection{Application: van der Waerden's theorem.}

We will now apply the theory that was presented in the previous subsection. The application will be a new proof of van der Waerden's theorem. This theorem says that, whenever the set $\mathbb N$ is colored with finitely many colors, there are arbitrarily long monochromatic arithmetic progressions.

\begin{theorem}[van der Waerden~\cite{vanderwaerden}]\label{vdw}
Given any coloring $c:\mathbb N\longrightarrow\{\r,\b\}$ and any finite number $l$, there are numbers $a,b\in\mathbb N$ such that the set $\{a,a+b,a+2b,\ldots,a+(l-1)b\}$ (which is an arithmetic progression of length $l$) is monochromatic.
\end{theorem}

\begin{proof}
Consider the set
\begin{equation*}
X=\{n+mx\mid n\in\mathbb N,m\in\mathbb N\cup\{0\}\},
\end{equation*}
the set of all polynomials (in the variable $x$) of degree either 0 or 1 with coefficients in $\mathbb N$. This set becomes a (commutative) semigroup when equipped with the usual addition operation. Furthermore, the set
\begin{equation*}
I=X\setminus\mathbb N=\{n+mx\mid n,m\in\mathbb N\}
\end{equation*}
of all polynomials in $X$ with degree {\it exactly} 1, is an ideal of the semigroup $X$. Thus, we have a situation where $X$ is a semigroup that can be written as the disjoint union $I\cup\mathbb N$, where $\mathbb N$ is a subsemigroup and $I$ is an ideal of $X$. All of these properties naturally carry over to the \v{C}ech--Stone compactification, and so we have that $\beta X$ is a right-topological semigroup that can be written as the disjoint union $\overline{I}\cup\overline{\mathbb N}$, where
\begin{equation*}
\overline{\mathbb N}=\{u\in\beta X\mid \mathbb N\in u\}
\end{equation*}
is a closed subsemigroup, and 
\begin{equation*}
\overline{I}=\{u\in\beta X\mid I\in u\}
\end{equation*}
is a closed ideal of $\beta X$.

Furthermore, for each $k\in\mathbb N\cup\{0\}$, we have an evaluation map
\begin{eqnarray*}
\ev_k:\ \ \ \ \ \ \ \ \ X & \longrightarrow & \mathbb N \\
n+mx & \longmapsto & n+mk
\end{eqnarray*}
which is a semigroup homomorphism and, when restricted to $\mathbb N$, is simply the identity map. Correspondingly, the Rudin--Keisler lifting of this mapping, $\ev_k:\beta X\longrightarrow\overline{\mathbb N}$ (which, as the reader might recall, maps $u$ to the ultrafilter $\{A\subseteq X\mid ev_k^{-1}[A]\in u\}$), will be a continuous semigroup homomorphism whose restriction to $\overline{\mathbb N}$ is simply the identity map.

With the previous two paragraphs setting up the tools, we now take our arbitrary coloring $c:\mathbb N\longrightarrow\{\r,\b\}$. Let $u\in\overline{\mathbb N}$ be a nonprincipal minimal idempotent (minimal within the closed subsemigroup $\overline{\mathbb N}$). Now let $v$ be a nonprincipal minimal idempotent (this time, minimal in all of $\beta X$) belonging to the closed ideal $\overline{I}$, with $v\leq u$. Fix an arbitrary $k\in\mathbb N\cup\{0\}$. Since $v\leq u$ and $\ev_k$ is a semigroup homomorphism, we have $\ev_k(v)\leq\ev_k(u)=u$. Now $\ev_k(v)\in\overline{\mathbb N}$ and $u$ was chosen to be minimal in $\overline{\mathbb N}$; therefore we must have $\ev_k(v)=u$.

Now, since $u$ is an ultrafilter and $\mathbb N\in u$, either $\{n\in\mathbb N\mid n\text{ is }\r\}\in u$ or $\{n\in\mathbb N\mid n\text{ is }\b\}\in u$. Assume without loss of generality that the first alternative holds, and we proceed to build an arithmetic progression of length $l$, all of whose elements are \tr. Let $A=\{n\in\mathbb N\mid n\text{ is }\r\}$. For each $k\in\{0,\ldots,l-1\}$, we have that $A\in u=\ev_k(v)=\{B\subseteq X\mid ev_k^{-1}[B]\in v\}$. In other words, $\ev_k^{-1}[A]\in v$. Since $v$ is an ultrafilter, and hence closed under finite intersections, we have
\begin{equation*}
I\cap\left(\bigcap_{k=0}^{l-1}\ev_k^{-1}[A]\right)\in v,
\end{equation*}
in particular this set is nonempty and so we can choose an element $a+bx\in I\cap\left(\bigcap_{k=0}^{l-1}\ev_k^{-1}[A]\right)$. Since $a+bx\in I$, we have $b\neq 0$. Furthermore, given a $k\in\{0,\ldots,l-1\}$, the fact that $a+bx\in\ev_k^{-1}[A]$ means that $a+kb=ev_k(a+bx)\in A$, in other words, $a+kb$ is \tr. This holds for every $k\in\{0,\ldots,l-1\}$ and so the whole set $\{a,a+b,\ldots,a+(l-1)b\}$ (which is an honest arithmetic progression of length $l$, since $b\neq 0$) consists only of \tr\ elements, and we are done.
\end{proof}

Suppose that we color $\mathbb N$ with colors \tr\ and \tb. By Theorem~\ref{vdw}, for each $l$ there is either an arithmetic progression of length $l$ that is completely \tr, or an arithmetic progression of length $l$ that is completely \tb. Say that $l$ is \tr\ or \tb\ according to whether the first or the second case holds. By the pigeonhole principle, either infinitely many $l$ are \tr, or infinitely many $l$ are \tb. But this means that either there are arbitrarily long arithmetic progressions that are \tr, or there are arbitrarily long arithmetic progressions that are \tb. Thus, Theorem~\ref{vdw} can be strengthened by saying that, whenever one partitions $\mathbb N$ into two pieces, then one of the pieces of the partition contains arbitrarily long arithmetic progressions.

Note that, for the three theorems exemplified in this article, the crucial steps in each of these proofs are when one partitions a set into two pieces and works with whichever piece belongs to a given ultrafilter. Since ultrafilters over a set $X$ contain a piece of every finite partition of $X$, it becomes clear that the same arguments go through if one considers colorings with any finite number of colors instead.

\begin{acknowledgment}{Acknowledgment.}
The author wishes to thank the two anonymous reviewers (whose careful reading and sound suggestions helped improve this paper nontrivially), and Fred Galvin for clarifying certain issues and providing some useful bibliographical pointers, especially regarding the history of the Galvin--Glazer argument. The author was supported by a postdoctoral fellowship from DGAPA--UNAM.
\end{acknowledgment}

\begin{biog}
\item[David Fern\'andez-Bret\'on] received his Ph.D. from York University in 2015. He currently is an assistant professor at Instituto Polit\'ecnico Nacional in Mexico City, after having held postdoctoral positions at the University of Michigan, the University of Vienna, Cinvestav, and Universidad Nacional Aut\'onoma de M\'exico
\begin{affil}
Instituto de Matem\'aticas, Universidad Nacional Aut\'onoma de M\'exico, \'Area de la Investigaci\'on Cient\'{\i}fica, Circuito Exterior, Ciudad Universitaria, Coyoac\'an, 04510, CDMX, Mexico \\
Current address: Escuela Superior de F\'{\i}sica y Matem\'aticas, Instituto Polit\'ecnico Nacional, Col. San Pedro Zacatenco, Alcald\'{\i}a Gustavo A. Madero, 07738, CDMX, Mexico. \\
{\tt dfernandezb@ipn.mx}
\end{affil}

\end{biog}
\vfill\eject

\end{document}